\documentclass{article}%
\usepackage{amsmath}
\usepackage{amsfonts}
\usepackage{amssymb}
\usepackage{graphicx}%
\addtocounter{footnote}{1}
\makeatletter
\@addtoreset{figure}{section}
\def\thefigure{\thesection.\@arabic\c@figure}
\def\fps@figure{h,t}
\@addtoreset{table}{bsection}
\def\thetable{\thesection.\@arabic\c@table}
\def\fps@table{h, t}
\@addtoreset{equation}{section}

\makeatother
\newtheorem{theorem}{Theorem}

\newtheorem{definition}[theorem]{Definition}
\newtheorem{example}[theorem]{Example}

\newtheorem{lemma}[theorem]{Lemma}

\newtheorem{remark}[theorem]{Remark}

\newenvironment{proof}[1][Proof]{\noindent\textbf{#1.} }{\ \rule{0.5em}{0.5em}}

\begin{document}

\title{\textbf{Antithetic variates in higher dimensions}}
\author{Sebastian del Ba\~{n}o Rollin$^{1}$ and Joan-Andreu L\'{a}zaro-Cam\'{\i}$^{2}%
\bigskip$\\$^{1}${\small Centre de Recerca Matem\`{a}tica. Apartat 50, E-08193 Bellaterra
}\\{\small Barcelona, Spain.}\\$^{2}${\small Imperial College London. Department of Mathematics.}\\{\small 180 Queen's Gate, SW7 2AZ, London, UK.}}
\date{}
\maketitle

\begin{abstract}
We introduce the concept of multidimensional antithetic as the absolute
minimum of the covariance function defined on the orthogonal group by
$A\mapsto\operatorname*{Cov}\left(  f(\xi),f(A\xi)\right)  $ where $\xi$ is a
standard $N$-dimensional normal random variable and $f:\mathbb{R}%
^{N}\rightarrow\mathbb{R}$ is an almost everywhere differentiable function.
The antithetic matrix is designed to optimise the calculation of
$\operatorname*{E}[f(\xi)]$ in a Monte Carlo simulation. We present an
iterative annealing algorithm that dynamically incorporates the estimation of
the antithetic matrix within the Monte Carlo calculation.

\end{abstract}

\makeatletter\addtocounter{footnote}{1} \footnotetext{e-mail:
\texttt{SebastiandB}{\texttt{@crm.cat}}} \addtocounter{footnote}{1}
\footnotetext{e-mail: \texttt{j.lazaro-cami}{\texttt{@imperial.ac.uk}}}
\makeatother

\bigskip


\bigskip

\section{Introduction}

The valuation of financial derivatives is based on the resolution of a
parabolic partial differential equation defined by the chosen dynamics for the
underlying assets subject to boundary conditions defined by the product (see
\cite{joshi}). These equations are rarely solvable explicitly and a numerical
method has to be chosen. The standard methods of choice in the industry are
resolution on grids and Monte Carlo. The applicability of the Monte Carlo
method is a consequence of the Feynman Kac theorem which solves a parabolic
PDE in terms of an expectation. In fact, in many of the more complex equity
products with a large dimensionality, the grid method is not efficient and
Monte Carlo is \emph{de facto} the only pricing method. In this approach the
price can be written as an expectation
\begin{equation}
m=\operatorname*{E}[f(\xi)], \label{antithetics 1}%
\end{equation}
where $\xi\sim\mathcal{N}(0,\operatorname*{Id}\nolimits_{N})$ is an
$N$-dimensional standard normal random variable describing a random path of
the underlying assets and $f:\mathbb{R}^{N}\rightarrow\mathbb{R}$ is a
measurable function representing the payoff of the derivative contract. The
Monte Carlo method is essentially a transcription of the strong law of large
numbers which claims that $m$ can be approximated by
\begin{equation}
\frac{1}{n}\sum_{i=1}^{n}f(\xi_{i}), \label{antithetics 2}%
\end{equation}
where $\xi_{1},\ldots,\xi_{n}$ are independent simulations of the random
variable $\xi$.

The main drawback of Monte Carlo methods is that they are usually
computationally demanding, often putting great strains on the capability of a
trading operation to properly monitor its risks and manage complex positions.
Therefore, it is of great value to design methods to improve the performance
of the Monte Carlo calculation. It is the aim of this article to present such
a method.

In order to speed up the Monte Carlo method we can seek to decrease the
simulation error. A measure of this error is given by $\sigma/\sqrt{n}$ for a
large enough number of simulations $n$ (\cite{Hammersley book}), where
$\sigma=\sqrt{\operatorname*{Var}\left[  f(\xi)\right]  }$ is the standard
deviation of $f\left(  \xi\right)  $. Like the expectation $m$, $\sigma$ is
usually unknown and needs to be estimated. {\bfseries\itshape Variance
reduction techniques} are methods that reduce this error by replacing
$f\left(  \xi\right)  $ by a different random variable which has the same
expectation but a smaller variance. Hopefully, this will ensure a faster
convergence of the Monte Carlo method. {\bfseries\itshape Antithetic variates
}is one such method.

Antithetic variates appear for the first time in the seminal work of
Hammersley, Morton, and Mauldon \cite{Hammersley,Hammersley2}. The crucial
idea of this procedure is to recycle the simulations of $\xi$ as samples of
$-\xi$, which has the same distribution as $\xi$. Therefore, we can
approximate $m$ in (\ref{antithetics 1}) by
\begin{equation}
\frac{1}{n}\sum_{i=1}^{n}\left(  \frac{f\left(  \xi_{i}\right)  +f\left(
-\xi_{i}\right)  }{2}\right)  . \label{antithetics 3}%
\end{equation}
The trivial equality
\[
\operatorname*{Var}\left[  f\left(  \xi\right)  +f\left(  -\xi\right)
\right]  =2\operatorname*{Var}\left[  f(\xi)\right]  +2\operatorname*{Cov}%
\left(  f(\xi),f(-\xi)\right)
\]
shows that if the dependence between $f(\xi)$ and $f(-\xi)$ is such that
$\operatorname*{Cov}\left(  f(\xi),f(-\xi)\right)  <0$, then the accuracy of
(\ref{antithetics 3}) will be greater than that of the crude Monte Carlo
(\ref{antithetics 2}).

Antithetic variates were developed in the case $\xi$ is one dimensional. For
higher dimensional normal variables, practitioners have often used antithetics
by changing the sign of some components of the random normal vector in a more
or less haphazard manner. In fact, there is a larger underlying group of
symmetries since, for any orthogonal matrix $A\in O(N)$, $A\xi$ is again a
standard normal vector. We can therefore extend the approach of Hammersley
\emph{et al.} to higher dimensions by replacing the role of $-\xi$ above by
$A\xi$ and approximating $m$ as
\begin{equation}
\frac{1}{n}\sum_{i=1}^{n}\frac{f\left(  \xi_{i}\right)  +f(A\xi_{i})}{2}.
\label{antithetics 26}%
\end{equation}
An antithetic method corresponds to the choice of $A\in O(N)$ and the optimal
antithetic is the matrix that minimizes the covariance function
\[
A\longmapsto\operatorname*{Cov}\left(  f(\xi),f(A\xi)\right)  .
\]
Note that, $O(N)$ being a compact group, this function will always have an
absolute minimum. The purpose of our work is to propose an algorithm to locate
this optimal antithetic matrix $A^{\ast}\in O(N)$, provided the covariance is
negative for some value $A\in O(N)$. The present paper is the first step of a
program whose aim is to describe optimal antithetics to price some popular
complex derivatives such as baskets, cliquets, Himalaya options and the like.

Finally, note that solving%
\[
\min_{A\in O(N)}\operatorname*{Cov}\left(  f(\xi),f(A\xi)\right)
\]
is equivalent to solving the simpler problem%
\begin{equation}
\min_{A\in O(N)}\operatorname*{E}\left[  f\left(  \xi\right)  f\left(
A\xi\right)  \right]  \label{antithetics 4}%
\end{equation}
and that both are well posed. Indeed, since $O(N)$ is compact, there exists a
$A^{\ast}\in O(N)$ such that $\min_{A\in O(N)}\operatorname*{E}\left[
f\left(  \xi\right)  f\left(  A\xi\right)  \right]  =\operatorname*{E}\left[
f\left(  \xi\right)  f\left(  A^{\ast}\xi\right)  \right]  $ provided that $f$
is continuous. More restrictively, we will assume throughout this paper that
$f$ is everywhere continuous and continuously differentiable except on a set
of zero (Lebesgue) measure.

\subsection{Contents}

The paper is structured as follows:

In Section \ref{section annealing algorithms} we recall some results about
simulated annealing algorithms. These algorithms are designed to find the
global minima of a given function defined on $\mathbb{R}^{r}$ by
stochastically perturbing a gradient based algorithm which, by itself, would
normally only converge to a \textit{local} minimum.

Section \ref{section Lie group methods} contains a detailed discussion on the
exponential covering map from the Lie algebra $\mathfrak{so}(N)$ to $O(N)$.
The coordinates induced by the exponetial map seem to be the best suited for
the optimization problem at hand.

In Section \ref{section antithetics revisited}, we introduce a novel iterative
simulated annealing algorithm adapted to the state space $O(N)$. This
algorithm provides a sequence $\{A_{k}\}_{k\in\mathbb{N}}\subset O(N)$ of
random variables converging in probability to the global minimum $A^{\ast}$ of
(\ref{antithetics 4}). The efficiency and performance of the algorithm are
checked in Section \ref{section examples}, where we use it to find the optimal
antithetic for an Asian call option and a covariance swap.

Finally, in Section \ref{section central limit}, we define a dynamical
antithetic technique where we use $\{A_{k}\}_{k\in\mathbb{N}}\subset O(N)$ to
estimate (\ref{antithetics 1}). More concretely, we define the sequence%
\begin{equation}
S_{n}:=\frac{1}{2n}\sum_{k=1}^{n}\left(  f\left(  \xi_{k}\right)  +f\left(
A_{k}\xi_{k}\right)  \right)  \label{antithetics 22}%
\end{equation}
and we prove that (\ref{antithetics 22}) converges almost surely to
$m=\operatorname*{E}[f(\xi)]$, and that $\sqrt{n}\left(  S_{n}-m\right)  $
converges in law to a normal variable $\mathcal{N}\left(  0,\sigma_{\ast}%
^{2}\right)  $ with variance%
\[
\sigma_{\ast}^{2}=\frac{1}{2}\left(  \operatorname*{Var}[f(\xi
)]+\operatorname*{Cov}\left(  f(\xi),f\left(  A^{\ast}\xi\right)  \right)
\right)  .
\]
That is, the estimated error of the Monte Carlo method (\ref{antithetics 22})
was reduced as much as it was possible using antithetics.

\subsection{Notation}

Throughout this article, $\left(  \Omega,\mathcal{F},P\right)  $ will denote
the underlying probability space, where $\mathcal{F}$ is a $\sigma$-algebra
and $P:\mathcal{F}\rightarrow\lbrack0,1]$ is a probability measure. The space
of all $\mathbb{R}^{N}$-valued random variables will be denoted by
$L_{\mathbb{R}^{N}}^{0}\left(  \Omega,P\right)  $. A standard Gaussian vector
$\xi:\Omega\rightarrow\mathbb{R}^{N}$ will be a random variable with a
probability density function given by
\[
\rho\left(  x\right)  =\frac{1}{(2\pi)^{N/2}}\operatorname*{e}%
\nolimits^{-\frac{1}{2}\left\Vert x\right\Vert ^{2}}%
\]
with respect to the Lebesgue measure $\lambda$. We will write $\xi
\sim\mathcal{N}(0,\operatorname*{Id}\nolimits_{N})$ to mean $\xi$ is a
standard Gaussian vector. On the other hand, $C_{\lambda}^{1}\left(
\mathbb{R}^{N}\right)  $ will denote the set of real functions $f:\mathbb{R}%
^{N}\rightarrow\mathbb{R}$ which are continuously differentiable except on a
set of zero Lebesgue measure.

\section{Simulated annealing algorithms\label{section annealing algorithms}}

In this section we review some of the existing methods to deal with the
problem of locating absolute minima for a function of the form
\[%
\begin{array}
[c]{rcl}%
\overline{U}:\mathbb{R}^{r} & \longrightarrow & \mathbb{R}\\
x & \longmapsto & \operatorname*{E}\left[  U\left(  x,\xi\right)  \right]  ,
\end{array}
\]
where $\xi\sim\mathcal{N}\left(  0,\operatorname*{Id}\nolimits_{N}\right)  $
is a Gaussian random vector and $U:\mathbb{R}^{r}\times\mathbb{R}%
^{N}\rightarrow\mathbb{R}$ is a measurable function continuously
differentiable with respect to the first $r$ entries. Actually, in our
original problem, $\overline{U}$ is defined on the compact Lie group $O(N)$
rather than $\mathbb{R}^{r}$; this additional feature will be dealt with later.

The absolute minima of the function $\overline{U}$ above will be zeroes of its
gradient field $\nabla\overline{U}$ which can, under regularity conditions, be
calculated as $\operatorname*{E}\left[  \nabla_{x}U(x,\xi)\right]  $.
Hopefully, one may expect to solve $\nabla\overline{U}\left(  x\right)  =0$ by
some numerical scheme such as the gradient method. However, the equation
\[
\nabla\overline{U}\left(  x\right)  =\operatorname*{E}\left[  \nabla
_{x}U\left(  x,\xi\right)  \right]  =0
\]
is defined through an expectation. This means that, in general, the gradient
$\nabla\overline{U}$ needs to be evaluated by a Monte Carlo simulation, which
can be too expensive. Recall that we want to solve (\ref{antithetics 4}) in
order to improve the efficiency of a (crude) Monte Carlo. Therefore,
calculating the gradient by Monte Carlo estimation in order to improve our
original (crude) Monte Carlo makes no sense. The solution is worse than the
problem. In order to overcome this difficulty, we try to work directly with
$\nabla_{x}U\left(  x,\xi\right)  $.

\subsection{Robbins-Monro algorithms}

Let $F\left(  \cdot,\xi\right)  :\mathbb{R}^{r}\rightarrow\mathbb{R}^{r}$ be a
measurable vector field depending on a Gaussian vector $\xi\in\mathcal{N}%
(0,\operatorname*{Id}_{N})$. {\bfseries\itshape Robbins-Monro algorithms}{
(\cite{Monro-Robbins})} are designed to solve an equation of the form
$\operatorname*{E}\left[  F(x,\xi)\right]  =0$ by means of a scheme algorithm
such as%
\begin{equation}
x_{n}=x_{n-1}+\gamma_{n}F\left(  x_{n-1},\xi_{n}\right)  ,
\label{antithetics 5}%
\end{equation}
where $\{\gamma_{n}\}_{n\in\mathbb{N}}$ is a non-negative sequence of real
numbers and $\{\xi_{n}\}_{n\in\mathbb{N}}$ are independent Gaussian vectors.
Observe that, in our particular case, $F=-\nabla_{x}U$. We will set
$\overline{F}(x)=\operatorname*{E}\left[  F(x,\xi)\right]  $.

\begin{theorem}
[{\cite[Theorem 1.4.26]{duflo}}]\label{Theorem Robbins-Monro}Assume that
$\overline{F}(x)$ has a zero $x^{\ast}$. Then the sequence defined by
(\ref{antithetics 5}) converges almost surely to $x^{\ast}$ for almost all
initial conditions $x_{0}$ provided that

\begin{description}
\item[\textbf{A1.}] $\left\langle x-x^{\ast},\overline{F}\left(  x\right)
\right\rangle <0$ for any $x\in\mathbb{R}^{r}$.

\item[\textbf{A2.}] $\sum_{n\in\mathbb{N}}\gamma_{n}=\infty$ and $\sum
_{n\in\mathbb{N}}\gamma_{n}^{2}<\infty$.

\item[\textbf{A3.}] $\operatorname*{E}\left[  \left\Vert F\left(  x_{n-1}%
,\xi_{n}\right)  \right\Vert ^{2}~|~\mathcal{F}_{n-1}\right]  \leq K\left(
1+\left\Vert \xi_{n}\right\Vert ^{2}\right)  $ a.s. for some constant $K>0$
where we set $\mathcal{F}_{n-1}$ to be the $\sigma$-algebra generated by the
random variables $\{\xi_{k}~|~k\leq n-1\}$.
\end{description}

Here $\left\langle \cdot,\cdot\right\rangle $ and $\Vert\cdot\Vert$ denote the
Euclidean scalar product and norm respectively.
\end{theorem}

\begin{remark}
An example of a sequence $\{\gamma_{n}\}_{n\in\mathbb{N}}$ verifying the
conditions of \textbf{A2} is $\gamma_{n}=\frac{c}{n}$ for some constant $c>0$.
\end{remark}

Among the hypotheses stated in Theorem \ref{Theorem Robbins-Monro},
\textbf{A1} is the most problematic. It means, on the one hand, that $x^{\ast
}$ is a unique zero of the vector field $\overline{F}$ and, on the other, that
$\left\Vert x-x^{\ast}\right\Vert $ is a Lyapunov function for $\overline{F}$
so that $x^{\ast}$ is an asymptotically stable equilibrium point. When
$\overline{F}=-\nabla\overline{U}$, this condition is sometimes hidden behind
the requirement that $\overline{U}\left(  x\right)  $ is convex and has a
unique minimum. Unfortunately, both conditions are too restrictive. In fact,
by Morse theory (see \cite{milnor}), the number of critical points of a
function defined on $O(N)$ such as our original covariance function has to be
at least the dimension of the total rational homology space of $O(N)$, which
is $2^{\left\lfloor N/2\right\rfloor +1}$ (see \cite[Corollary III.3.15]%
{mimura}).

There are variants of the Robbins-Monro algorithms that are guaranteed to
converge to critical points of $\overline{F}$ which are actual minima of
$\overline{U}$. However, there is no certainty that these will be absolute
minima. In the next subsection, we will discuss simulated annealing methods
which will converge to global minima.

Robbins-Monro algorithms combined with {\bfseries\itshape importance sampling}
methods have also been successfully used in the context of derivative pricing
to reduce the variance of Monte Carlo simulations (see \cite{arouna 1},
\cite{arouna 2}, \cite{dufresne}, and \cite{fu-su}). For example, if the price
of a derivative is $\operatorname*{E}[f(\xi)]$ with $\xi\sim\mathcal{N}%
(0,\operatorname*{Id}\nolimits_{N})$ as in (\ref{antithetics 1}), then a
simple change of variables shows that
\[
\operatorname*{E}[f(\xi)]=\operatorname*{E}\left[  f(\xi+x)\mathrm{e}%
^{-x\cdot\xi-\frac{1}{2}\Vert x\Vert^{2}}\right]  .
\]
However, the variance of the variable under the second expectation operator
now depends on $x$. So we will achieve faster convergence by choosing $x$ that
minimises $\operatorname*{Var}[f(\xi+x)\mathrm{e}^{-x\cdot\xi-\frac{1}{2}\Vert
x\Vert^{2}}]$ or, equivalently,
\begin{equation}
\overline{U}(x)=\operatorname*{E}[f^{2}\left(  \xi\right)  \operatorname*{e}%
\nolimits^{-x\cdot\xi+\frac{1}{2}\left\Vert x\right\Vert ^{2}}].
\label{importance sampling}%
\end{equation}
If turns out this function is convex and thus has a unique minimum. The
Robbins-Monro algorithm (suitable modified with a truncation method to ensure
condition \textbf{A3}) can then be applied to $F(x,\xi)=\nabla_{x}[f^{2}%
(\xi)\mathrm{e}^{-x\cdot\xi+\frac{1}{2}\Vert x\Vert^{2}}]$ to yield the
optimal value of $x\in\mathbb{R}^{r}$ (\cite{arouna 1}, \cite{arouna 2}). Note
that in our problem the function to be minimised is not a measure change, and
that it does not have the growth properties at infinity as a function of $x$
of (\ref{importance sampling}).

\subsection{Simulated annealing algorithms}

As we have mentioned, our minimisation problem will give rise to multiple
critical points and thus is not suited to the Robbins-Monro scheme. A method
that has been devised to deal with multiple minima is Simulated Annealing.
This is a technique inspired in the metallurgy where a metal alloy is heated
to pull it out of a equilibrium state, a local minimum of the energy, and then
slowly cooled to allow atoms to diffuse into the lowest energy state where the
system has some optimal physical property. The analogue of this heat injection
in the Robbins-Monro algorithm (\ref{antithetics 5}) is the addition of an
extra source of randomness that \textit{hits} the approximating sequence out
of local minima. In their most general form, {\bfseries\itshape simulated
annealing algorithms} are written as
\begin{equation}
x_{n+1}=x_{n-1}+\gamma_{n}F\left(  x_{n-1},\xi_{n}\right)  +b_{n}\zeta_{n},
\label{antithetics 6}%
\end{equation}
where $\{b_{n}\}_{n\in\mathbb{N}}\subset\mathbb{R}$ is a real sequence, the
annealing temperature scheme, and $\{\zeta_{n}\}_{n\in\mathbb{N}}$ a second
sequence of i.i.d Gaussian random vectors independent from the $\{\xi
_{n}\}_{n\in\mathbb{N}}$. In order to study their convergence, first of all,
we need to introduce some notation.

Let $B_{x,y}$ denote the set of continuous paths $\varphi:[0,1]\rightarrow
\mathbb{R}^{r}$ starting at $x\in\mathbb{R}^{r}$ and ending at $y\in
\mathbb{R}^{r}$. Let $U(x,\xi)$ be as above and set $F(x,\xi)=-\nabla
_{x}U(x,\xi)$ and $\overline{U}(x)=\operatorname*{E}[U(x,\xi)]$. Recall that
our intention is to locate absolute minima of $\overline{U}$. Let
$\underline{S}$ be the set where the absolute minimum is achieved
\[
\underline{S}=\{x\in\mathbb{R}^{r}~|~\overline{U}\left(  x\right)  =\min
_{y\in\mathbb{R}^{r}}\overline{U}\left(  y\right)  \}
\]
and $S=\{x\in\mathbb{R}^{r}~|~\nabla\overline{U}\left(  x\right)  =0\}$ the
set of critical points of $\overline{U}$ which we will assume has finitely
many connected components. For any $l>0$, let $B\left(  \underline
{S},l\right)  $ be the set of points which are at a distance less than $l$
from $\underline{S}$. Define%
\begin{equation}
I\left(  x,y\right)  =\inf_{\varphi\in B_{x,y}}2\left(  \max_{0\leq t\leq
1}\overline{U}\left(  \varphi(t)\right)  -\overline{U}(x)\right)  \text{ and
}d^{\ast}=\max_{x\in S\backslash\underline{S}}\min_{y\in\underline{S}}I\left(
x,y\right)  . \label{antithetics 18}%
\end{equation}
Observe that $d^{\ast}$ is a measure of how \textit{oscillatory} the objective
function $\overline{U}$ is.

The next theorem provides sufficient conditions for guaranteeing the
convergence of simulated annealing algorithms. It is extracted from
\cite[Theorem 5.2]{fang-gong-quian}.

\begin{theorem}
\label{Theorem Fang-Gong-Quian}Let $\{r_{n}\}_{n\in\mathbb{N}}\subset
\mathbb{R}$ and $\{h_{n}\}_{n\in\mathbb{N}}\subset\mathbb{R}$ be the sequences
of real numbers defined by $r_{n}=1/n^{\gamma}$, $0<\gamma<1$, and
$h_{n}=d/((1-\gamma)\ln n)$, where $d>d^{\ast}$. For $\{\xi_{n}\}_{n\in
\mathbb{N}}$, $\{\zeta_{n}\}_{n\in\mathbb{N}}$ two independent sequences of
i.i.d Gaussian random vectors define an iteration scheme by
\begin{equation}
x_{n}=x_{n-1}-r_{n}\nabla_{x}U\left(  x_{n-1},\xi_{n}\right)  +\sqrt
{r_{n}h_{n}}\zeta_{n}. \label{antithetics 8}%
\end{equation}
Assume the function $\overline{F}(x)=-\operatorname*{E}[\nabla_{x}U\left(
x,\xi\right)  ]$ satisfies the following properties:

\begin{description}
\item[B1.] $\lim\sup_{\left\Vert x\right\Vert \rightarrow\infty}%
\frac{\left\Vert \overline{F}\left(  x\right)  \right\Vert }{\left\Vert
x\right\Vert }\leq M_{1}<\infty$,

\item[B2.] $\lim\sup_{\left\Vert x\right\Vert \rightarrow\infty}%
\frac{\left\langle \overline{F}\left(  x\right)  ,x\right\rangle }{\left\Vert
x\right\Vert ^{2}}\leq-c<0$, and

\item[B3.] $\lim\sup_{\left\Vert x\right\Vert \rightarrow\infty}%
\frac{\operatorname*{E}\left[  \left(  F\left(  x,\xi\right)  -\overline
{F}\left(  x\right)  \right)  ^{2}\right]  }{\left\Vert x\right\Vert ^{2}}\leq
M_{2}<\infty$
\end{description}

\noindent for some positive constants $M_{1}$, $M_{2}$, and $c$. Then, for any
$l>0$,
\begin{equation}
P\left(  \{x_{n}\in B\left(  \underline{S},l\right)  \}\right)  \rightarrow
1\text{ as }n\rightarrow\infty\label{antithetics 16}%
\end{equation}
uniformly for any initial condition $x_{0}$ in an arbitrary compact set.
\end{theorem}

\begin{remark}
\normalfont \ \ 

\begin{enumerate}
\item If the absolute minimum $x^{\ast}$ of $\overline{U}$ is unique, then
(\ref{antithetics 16}) implies that the sequence $\{x_{n}\}_{n\in\mathbb{N}}$
converge in probability to $x^{\ast}$.

\item In fact, Fang \emph{et al.} have proved a more general version of
Theorem \ref{Theorem Fang-Gong-Quian} for a wider spectrum of sequences
$\{r_{n}\}_{_{n\in\mathbb{N}}}$ and $\{h_{n}\}_{n\in\mathbb{N}}$ (see
\cite{fang-gong-quian} for details). For example, they prove that Theorem
\ref{Theorem Fang-Gong-Quian} holds when $r_{n}=b/n$ and $h_{n}=d/\ln(\ln n)$,
$b>0$, $d>0$, and $n\in\mathbb{N}$. Simulated annealing algorithms with
coefficients $-b/n$ and $\sqrt{bd/\left(  n\ln(\ln n)\right)  }$ have already
been considered in \cite{Gelfand-Mitter,Gelfand-Mitter-metropolys} in the
particular case $\nabla_{x}U\left(  x_{n-1},\xi_{n}\right)  =\nabla
_{x}V\left(  x_{n-1}\right)  +\xi_{n}$ with $V:\mathbb{R}^{r}\rightarrow
\mathbb{R}$ a deterministic function.

\item Hypotheses \textbf{B1}, \textbf{B2}, and \textbf{B3} are weaker than the
corresponding hypotheses \textbf{A1}, \textbf{A2}, and \textbf{A3} in the
Robbins-Monro algorithm. This is because the former only impose restrictions
\textit{at infinity}. For example, if we know that the absolute minimum of the
function $\overline{U}$ lies in a bounded domain, we can modify $\overline{U}$
far away from that domain so that \textbf{B1}, \textbf{B2}, and \textbf{B3}
automatically hold and $\{x_{n}\}_{n\in\mathbb{N}}$ in (\ref{antithetics 8})
converges. If we are certain that the minimum is within the ball $B\left(
x^{\ast},R\right)  $ of radius $R$ centered at some $x^{\ast}\in\mathbb{R}%
^{r}$, then we can add a penalty function appropriately on a narrow spherical
shell around $\left\Vert x\right\Vert =R$ such that $-\nabla_{x}\overline{U}$
strongly points towards $B\left(  x^{\ast},R\right)  $.

\item As far as the rate of convergence of (\ref{antithetics 8}) is concerned,
it is also proved in \cite{fang-gong-quian} that there exists a constant
$\delta>0$ such that, for any $\varepsilon>0$,
\[
P\left(  \{x_{n}\in B\left(  \underline{S},l\right)  \}\right)  \geq
1-\varepsilon\text{ \ if }n>\exp\left(  -\frac{d}{(1-\gamma)\delta}%
\ln(2\varepsilon)\right)  .
\]

\end{enumerate}
\end{remark}

\subsection{Continuous simulated annealing}

Simulated annealing is often presented in the literature in its continuous
version. The stochastically perturbed iterative algorithm is then replaced by
a stochastic differential equation with drift $-\nabla_{x}\overline{U}$ the
negative of the gradient of the function $\overline{U}\in C^{1}\left(
\mathbb{R}^{r}\right)  $ we wish to minimise, and diffusion coefficient
decaying to zero with time as in the annealing method. Unfortunately, as it
will become clearer later on, this procedure will only work for deterministic
functions; that is, functions which do not depend on a random variable as is
indeed our case. However, we wish to review time-continuous simulated
annealing for the benefit of a more complete exposition. Most of the results
quoted here are extracted from \cite{josef-fabrice} where Baudoin, Hairer, and
Teichmann study the Ornstein-Uhlenbeck process on a compact Lie group and its
properties to design efficient simulated annealing schemes. This recent paper
improves considerably the efficiency of simulated annealing techniques
developed so far (see \cite{Holley}). The reader is also encouraged to check
with \cite[Chapter 5]{bayer PhD thesis} for a more comprehensive approach.

Let $G$ be a compact Lie group and let $\mathcal{L}=\frac{1}{2}\sum_{i=1}%
^{d}V_{i}\circ V_{i}$ be a second order differential operator acting on
$L^{2}\left(  G,\mu_{G}\right)  $, where $\mu_{G}$ is the
%(right-invariant)
Haar measure and, for any $i=1,...,d$, $V_{i}\in\mathfrak{g}$ is a left
invariant vector field. We assume that H\"{o}rmander's hypoelliptic condition
holds, i.e., that the Lie sub-algebra generated by $\{V_{1},...,V_{d}\}$
coincides with $\mathfrak{g}$. In this context, the {\bfseries\itshape
carr\'{e} du champ operator} $\Gamma$ is defined as
\[
\Gamma\left(  g,f\right)  =\mathcal{L}\left(  fg\right)  -f\mathcal{L}\left(
g\right)  -g\mathcal{L}\left(  f\right)  ,~~f,g\in L^{2}\left(  G,\mu
_{G}\right)  .
\]

Let $\overline{U}\in C^{1}\left(  G\right)  $ be a differentiable function and
let $\varepsilon\in(0,1]$ be the annealing temperature parameter. Assume that
the following integral is finite
\[
Z_{\varepsilon}:=\int_{G}\operatorname*{e}\nolimits^{-\overline{U}%
(g)/\varepsilon^{2}}d\mu_{G}(g)<\infty.
\]
Then we can define the {\bfseries\itshape Gibbs measure} $\mu_{\varepsilon}$
by $\mu_{\varepsilon}\left(  B\right)  :=\frac{1}{Z_{\varepsilon}}\int
_{B}\operatorname*{e}\nolimits^{-\overline{U}(g)/\varepsilon^{2}}d\mu_{G}(g)$,
$B\in\mathcal{B}\left(  G\right)  $. Intuitively this measure concentrates on
the minima of $\overline{U}$ as the temperature $\varepsilon$ falls to zero.
The Gibbs measure is invariant under the differential operator $\mathcal{L}%
_{\varepsilon}=\varepsilon^{2}\mathcal{L}-\frac{1}{2}\Gamma\left(
\overline{U},\cdot\right)  $, which is the {\bfseries\itshape infinitesimal
generator} of the stochastic differential equation
\begin{equation}
dX_{t}^{g}=V_{0}\left(  X_{t}^{g}\right)  dt+\varepsilon\sum_{i=1}^{d}%
V_{i}\left(  X_{t}^{g}\right)  \delta B_{t}^{i} \label{antithetics 14}%
\end{equation}
where $V_{0}=-\frac{1}{2}\Gamma\left(  \overline{U},\cdot\right)
\in\mathfrak{X}\left(  G\right)  $ represents the negative gradient vector of
the function $\overline{U}$. This is the continuous version of the simulated
annealing process. Note that if the temperature $\varepsilon$ is set to zero,
we are left with an ordinary differential equation where the \textit{carr\'{e}
du champ} $\Gamma(\overline{U},\cdot)$ plays the r\^{o}le of the ordinary
gradient. Its flow, hopefully, will evolve to the closest minimum of
$\overline{U}$. The addition of appropriately selected annealing schedule
function $\varepsilon(t)$ will ensure we move to an absolute minimum.

In \cite{josef-fabrice} the following remarkable result is proved:

\begin{theorem}
Assume the annealing heat function is given by%
\[
\varepsilon\left(  t\right)  =\frac{c}{\sqrt{\ln\left(  R+t\right)  }}%
\]
for positive constants $c,R>0$. Then, under mild conditions on $\overline
{U}\in C^{1}\left(  G\right)  $,
\begin{equation}
P\left(  \left\{  X_{t}^{g}\in B_{\delta}\right\}  \right)  \leq M\sqrt
{\mu_{\varepsilon(t)}\left(  B_{\delta}\right)  }, \label{antithetics 15}%
\end{equation}
where $M>0$, $B_{\delta}:=\left\{  g\in G~|~\overline{U}(g)>\overline{U}%
_{0}+\delta\right\}  $, and $\overline{U}_{0}$ is the absolute minimum of
$\overline{U}$. In particular, provided that there exists only one element
$g_{0}\in G$ such that $\overline{U}(g_{0})=\overline{U}_{0}$,
(\ref{antithetics 15}) implies that%
\[
\lim_{t\rightarrow\infty}\operatorname*{E}\left[  f\left(  X_{t}^{g}\right)
\right]  =f\left(  g_{0}\right)
\]
for any continuous bounded test function $f\in C\left(  G\right)  $.
\end{theorem}

In other words, in order to find the global minima of $\overline{U}$ we can
simulate numerically the stochastic differential equation
(\ref{antithetics 14}) and approximate its invariant measure $\mu
_{\varepsilon}$ which, as time goes by, concentrates around $g_{0}\in G$.

Unfortunately, our original problem does not fit directly in this framework.
The drift vector field $V_{0}=\frac{1}{2}\Gamma\left(  \overline{U}%
,\cdot\right)  $ in (\ref{antithetics 14}) is assumed to be deterministic, it
cannot depend on an independent Gaussian noise, as happens in our case:
according to (\ref{antithetics 4}), $U=f\left(  \xi\right)  f\left(
A\xi\right)  $, where $\xi\sim\mathcal{N}(0,\operatorname*{Id}\nolimits_{N})$,
$A\in O(N)$, and $f\in C_{\lambda}^{1}(\mathbb{R}^{N})$. Removing this
stochastic dependence would imply working with $\overline{U}%
(A)=\operatorname*{E}[f\left(  \xi\right)  f\left(  A\xi\right)  ]$. As we
already argued, it is computationally inconvenient to take the expectation at
this point.

\section{Lie group methods\label{section Lie group methods}}

We have so far reviewed the simulated annealing algorithms available to
globally minimize functions defined on $\mathbb{R}^{r}$ through an
expectation. However, our optimization problem (\ref{antithetics 4}) does not
take place on an Euclidean space but on a Lie group, namely, the orthogonal
group $O(N)$. Therefore, we need to adapt the results of Section
\ref{section annealing algorithms} to $O(N)$ by taking suitable coordinates.
As it is customary on Lie groups, we will take local coordinates by means of a
local diffeomorphism from the Lie algebra, $\varphi:\mathfrak{o}\left(
N\right)  \rightarrow O(N)$, which is a Euclidean space. In this section we
are going consider two different choices of $\varphi$: the Cayley transform
and the exponential map. These are two of the most used coordinate patches in
the design of numerical integrators for ordinary differential equations on Lie
groups (\cite{lie group methods}). This section aims at giving explicit
expressions for the gradient of a smooth function $F:O(N)\rightarrow
\mathbb{R}$ when composed with the local coordinates $\varphi$. Such a
gradient will be employed later in our simulated annealing algorithm.

Let $G$ be a Lie group and $\mathfrak{g}$ its Lie algebra. The tangent bundle
$\tau_{G}:TG\rightarrow G$ is trivial meaning that it is isomorphic to the
product $G\times\mathfrak{g}$. The identification $TG=G\times\mathfrak{g}$ can
be carried out by means of an isomorphism $\rho:TG\rightarrow G\times
\mathfrak{g}$ induced by right translations. If $R_{g}:G\rightarrow G$ denotes
the map defined for $g\in G$ by $R_{g}(h)=hg$, then $\rho\left(  v\right)
=(g,T_{g}R_{g^{-1}}(v))$, where $g=\tau_{G}\left(  v\right)  $. We refer to
this trivialisation as the {\bfseries\itshape space coordinates} on the
tangent bundle.

Given a smooth function $F:G\rightarrow\mathbb{R}$, the tangent map
$TF:TG\rightarrow T\mathbb{R}$ it induces will be noted $TF^{sc}$ when written
in space coordinates. That is,
\[%
\begin{array}
[c]{rcl}%
TF^{sc}:G\times\mathfrak{g} & \longrightarrow & \mathbb{R}^{2}\\
\left(  g,X\right)  & \longmapsto & \left(  F\left(  g\right)  ,T_{g}F\circ
T_{e}R_{g}(X)\right)  .
\end{array}
\]
Equivalently, if $\varphi:\mathfrak{g}\rightarrow G$ is a local diffeomorphism
from a neighborhood of $0\in\mathfrak{g}$, $T\varphi^{sc}:\mathfrak{g}%
\times\mathfrak{g}\rightarrow G\times\mathfrak{g}$ stands for the tangent map
$T\varphi:\mathfrak{g}\times\mathfrak{g}\rightarrow TG$ in space coordinates,
i.e., $T_{X}\varphi^{sc}=T_{\varphi(X)}R_{\varphi(X)^{-1}}\circ T_{X}\varphi$,
$X\in\mathfrak{g}$. It can be immediately checked that%
\[
T\left(  F\circ\varphi\right)  =TF^{sc}\circ T\varphi^{sc}.
\]

In order to minimise a function on $G$ or, via a coordinate chart, on
$\mathfrak{g}$, it will be useful to have a notion of gradient field. Assume
therefore that we are given an arbitrary metric $\left\langle \cdot
,\cdot\right\rangle $ on $\mathfrak{g}$. We will describe the gradient of a
function $F:G\rightarrow\mathbb{R}$ in space coordinates. For $X,Y\in
\mathfrak{g}$, the gradient of $F\circ\varphi:\mathfrak{g}\rightarrow
\mathbb{R}$ satisfies%
\begin{align*}
\left\langle \nabla\left(  F\circ\varphi\right)  \left(  X\right)
,Y\right\rangle  &  =d\left(  F\circ\varphi\right)  \left(  X\right)  \left(
Y\right)  =pr_{2}\circ T_{X}\left(  F\circ\varphi\right)  \left(  Y\right) \\
&  =pr_{2}\circ T_{\varphi(X)}F^{sc}\circ T_{X}\varphi^{sc}(Y)
\end{align*}
where $pr_{2}:\mathbb{R}\times\mathbb{R}\rightarrow\mathbb{R}$ denotes the
projection onto the second factor. Therefore, if $\{Y_{i}\}_{i=1,...,\dim
(\mathfrak{g})}$ is an orthonormal basis, then%
\begin{equation}
\nabla\left(  F\circ\varphi\right)  \left(  X\right)  =\sum_{i=1}%
^{\dim(\mathfrak{g})}\left(  pr_{2}\circ T_{\varphi(X)}F^{sc}\circ
T_{X}\varphi^{sc}(Y_{i})\right)  Y_{i}. \label{antithetics 11}%
\end{equation}

\begin{example}
[Canonical coordinates of the first kind]\normalfont Let $\varphi
:\mathfrak{g}\rightarrow G$ be $\varphi\left(  X\right)  :=\exp\left(
X\right)  g$ for some $g\in G$. The exponential map of a Lie group is a local
diffeomorphism from a neighborhood of $0\in\mathfrak{g}$ onto a neighborhood
of $g\in G$. Using these coordinates, it can be checked that%
\begin{equation}
T_{X}\varphi^{sc}=\frac{\exp\left(  \operatorname*{ad}_{X}\right)
-\operatorname*{Id}}{\operatorname*{ad}_{X}}=\sum_{j\geq0}\frac{1}%
{(j+1)!}\operatorname*{ad}\nolimits_{X}\circ\overset{j)}{\ldots}%
\circ\operatorname*{ad}\nolimits_{X} \label{antithetics 12}%
\end{equation}
(see \cite{lie group methods,Engo}). Canonical coordinates of the first kind
are convenient for nilpotent Lie algebras because then the series
(\ref{antithetics 12}) becomes a finite sum. Unfortunately, $\mathfrak{so}%
\left(  N\right)  $ is not nilpotent. However, these coordinates are useful
for $SO(3)$ because, in this particular case, the exponential can be easily
computed by means of the \textit{Rodrigues formula}. Indeed, if%
\[
X=%
\begin{pmatrix}
0 & -a & b\\
a & 0 & -c\\
-b & c & 0
\end{pmatrix}
\in\mathfrak{so}\left(  3\right)  ,
\]
one can prove that%
\[
\exp\left(  X\right)  =\operatorname*{Id}+\frac{\sin(\sigma)}{\sigma}%
X+\frac{1-\cos(\sigma)}{\sigma^{2}}X^{2},
\]
where $\sigma=\sqrt{a^{2}+b^{2}+c^{2}}$, and%
\[
T_{X}\exp=\operatorname*{Id}+\frac{1-\cos(\sigma)}{\sigma}X+\frac{\sigma
-\sin(\sigma)}{\sigma^{3}}X^{2}.
\]
For $n>3$ several methods have been devised to calculate the exponential map
by other means than truncation of the series $\sum_{n\geq\infty}\frac{1}%
{n!}X^{n}$, which often leads to numerical error (see \cite{buono-exponential}%
,\cite{gallier-exponential}).
\end{example}

\section{Higher dimensional antithetic
variates\label{section antithetics revisited}}

We now present an annealing type method to find the optimal antithetic matrix
$A^{\ast}\in O(N)$ as defined in the introduction. Recall that an optimal
antithetic is an absolute minimum of the function $A\mapsto\operatorname*{Cov}
\left(  f\left(  \xi\right)  ,f\left(  A\xi\right)  \right)  $, where $f\in
C_{\lambda}^{1}(\mathbb{R}^{N})$ is the payout function depending on a path
described by a normal vector $\xi\sim\mathcal{N}\left(  0,\operatorname*{Id}%
\nolimits_{N}\right)  $. This is equivalent to minimising $\overline
U(A)=\operatorname*{E}\left[  f\left(  \xi\right)  f\left(  A\xi\right)
\right]  $. In the notation of previous sections we write $U(A,\xi)=f\left(
\xi\right)  f\left(  A\xi\right)  $.

The plan is to coordinate $SO(N)$ by means of the exponential map $\exp(Y)$,
$Y\in\mathfrak{so}(N)$, introduce an iterative algorithm, and study its
convergence in the light of Theorem \ref{Theorem Fang-Gong-Quian}.
Unfortunately, $U(\exp(Y),\xi)$ have not the appropriate behaviour at infinity
as specified by hypothesis \textbf{B2} in Theorem
\ref{Theorem Fang-Gong-Quian}. Indeed, let $\left\langle X,Y\right\rangle
:=\frac{1}{2}\operatorname*{trace}(XY^{\top})$ be the standard Euclidean
product on $\mathfrak{so}(N)$ and let $\left\Vert \cdot\right\Vert $ its
associated norm, $Y,X\in\mathfrak{so}(N)$. If $\overline{F}=-\nabla
\operatorname*{E}[f\left(  \xi\right)  f\left(  \exp(\cdot)\xi\right)  ]$
denotes the gradient of $\overline{U}\circ\exp$ with respect to an orthonormal
basis of $\mathfrak{so}(N)$ then, since $\overline{U}$ is bounded,%
\[
\lim_{\Vert Y\Vert\rightarrow\infty}\frac{\langle\overline{F}(Y),Y\rangle
}{\left\Vert Y\right\Vert ^{2}}=0,
\]
which is not strictly negative as \textbf{B2} requires. In order that our
algorithm satisfies this and the rest of hypotheses, we will therefore modify
$\overline{U}(\exp(Y))$ far away from a bounded set with a penalty function.

Before stating the main results of this section, we need an auxiliary lemma.
Recall that the {\bfseries\itshape rank} of a Lie group is the dimension of a
maximal torus and that $SO(N)$ has rank $\left\lfloor \left.  N\right/
2\right\rfloor $, where $\left\lfloor \cdot\right\rfloor $ stands for the
integer part of a real number. We define $R:=\left\lfloor \left.  N\right/
2\right\rfloor $ for the sake of a simpler notation.

\begin{lemma}
\label{SurjectingBall} The compact ball $B(0,\pi\sqrt{R})\subset
\mathfrak{so}(N)$ surjects onto $SO(N)$ by the exponential map.
\end{lemma}

\begin{proof}
It is a well know fact that any matrix in $SO(N)$ is conjugated to an element
in the maximal torus of matrices with $\left\lfloor N/2\right\rfloor $
diagonal blocks of the form%
\[%
\begin{pmatrix}
\cos\theta & -\sin\theta\\
\sin\theta & \cos\theta
\end{pmatrix}
,~\theta\in\lbrack-\pi,\pi],
\]
(see \cite[Chapter IV Theorem 1.6]{tom Dieck}). It is immediate to see that
these rotations can be written as $\exp%
\begin{pmatrix}
0 & -\theta\\
\theta & 0
\end{pmatrix}
$. Since conjugation commutes with exponentiation, we conclude that, if%
\[
K:=\left\{  \operatorname*{diag}\left(  S_{1},...,S_{R}\right)  \in
\mathfrak{so}(N)~|~S_{i}=%
\begin{pmatrix}
0 & -\theta_{i}\\
\theta_{i} & 0
\end{pmatrix}
,~\theta_{i}\in\lbrack-\pi,\pi],~i=1,...,R\right\}  ,
\]
then the compact set $\left\{  AkA^{\top}~|~A\in SO(n),~k\in K\right\}
\subset\mathfrak{so}(N)$ surjects onto $SO(n)$ by the exponential map. The
lemma follows because, for any $A\in SO(n)$ and any $k\in K$, we have
\[
\Vert AkA^{\top}\Vert^{2}=\frac{1}{2}\operatorname*{trace}\left(  Akk^{\top
}A^{\top}\right)  =\frac{1}{2}\operatorname*{trace}(kk^{\top})\leq R\pi^{2}%
\]

\end{proof}

As stated above, we will use a penalty function to ensure that our algorithm
converges. We use the notation $\mathbf{1}_{B}(Y)$ for the characteristic
function of a set $B\subseteq\mathfrak{so}\left(  N\right)  $.

\begin{definition}
Let $P(Y)$ be the function defined on $\mathfrak{so}(N)$ by
\[
P(Y)=\Vert Y\Vert^{2}\cdot\mathbf{1}_{B\left(  0,\pi\sqrt{R}\right)  ^{c}}(Y)
\]

\end{definition}

The algorithm we present in the following theorem requires some notation. Let
$E_{ij}\in\mathfrak{so}(N)$ be the matrix with $+1$ in the $(i,j)$th position,
$-1$ in the $(j,i)$th position and zero elsewhere. The set $\left\{
E_{ij}\right\}  _{0<i<j\leq N}$ is a basis of $\mathfrak{so}(N)$, orthogonal
with respect to the aforementioned standard scalar product.

\begin{theorem}
Let $\left\{  \xi_{n}\right\}  _{n\in\mathbb{N}}$ and $\left\{  \zeta
_{n}\right\}  _{n\in\mathbb{N}}$ be sequences of independent standard normal
$N$-dimensional vectors, $A_{0}\in O(N)$ and let $\varphi$ be the exponential
function centered at $A_{0}$, $\varphi(Y)=\exp(Y)A_{0}$. Set $Y_{0}=0$ and
$Z_{0}=0$. Then the sequence $\{A_{n}\}_{n\in\mathbb{N}}$ defined by
\begin{subequations}
\label{algorithm steps}%
\begin{align}
Z_{n}  &  =\sum_{i<j}\left(  f(\xi_{n})\nabla f^{\intercal}|_{A_{n-1}\xi_{n}%
}T_{Y_{n-1}}\varphi^{sc}(E_{ij})A_{n-1}\xi_{n}\right)  E_{ij}\in
\mathfrak{so}(N)\label{algorithm step 1}\\
Y_{n}  &  =Y_{n-1}-\frac{1}{n^{\gamma}}\left(  Z_{n}+2Y_{n-1}\cdot
\mathbf{1}_{B\left(  0,\pi\sqrt{R}\right)  ^{c}}\left(  Y_{n-1}\right)
\right)  +\sqrt{\frac{d}{n^{\gamma}(1-\gamma)\ln n}}\zeta_{n}%
\label{algorithm step 2}\\
A_{n}  &  =\exp(Y_{n})A_{0} \label{algorithm step 3}%
\end{align}
converges in probability to the optimal antithetic $A^{\ast}$ (or the
antithetic locus if there are multiple minima) in the connected component of
$SO(N)$ containing $A_{0}$. Here $d>0$ is the same as in Theorem
\ref{Theorem Fang-Gong-Quian}.
\end{subequations}
\end{theorem}

\begin{proof}
The proof consists in showing that the iterative process $\{Y_{n}%
\}_{n\in\mathbb{N}}$ is the simulated annealing algorithm described in Theorem
\ref{Theorem Fang-Gong-Quian} applied to the function $U(\xi,\exp(Y))+P(Y)$,
$Y\in\mathfrak{so}\left(  N\right)  $.

Since $U(\xi,\exp(Y))+P(Y)$ coincides with $U(\xi,\exp(Y))$ on $B(0,\pi
\sqrt{R})$ and is strictly greater outside this ball, the minima of these
functions will coincide in $B(0,\pi\sqrt{R})$. Thanks to Lemma
\ref{SurjectingBall}, the exponential of the absolute minima of $U(\xi
,\exp(Y))+P(Y)$ yield the optimal antithetics.

The Robbins-Monro gradient term in the annealing procedure is
\begin{equation}
F(\xi_{n},Y_{n-1})=-\nabla_{Y}\left(  U(\xi_{n},\exp(Y))+P(Y)\right)
|_{Y=Y_{n-1}}. \label{RMGradient}%
\end{equation}
In order to calculate it, let $\varphi(Y)=\exp(Y)A_{0}$ and write $U_{\xi
}(A):=U(\xi,A)$, $A\in O(n)$. The tangent map $T_{A}U_{\xi}:T_{A}%
O(N)\rightarrow T_{U_{\xi}(A)}\mathbb{R}\simeq\mathbb{R}$ can be described in
space coordinates on $X\in\mathfrak{so}(N)$ by using the chain rule in the
following manner
\begin{align*}
T_{A}U_{\xi}^{sc}\left(  X\right)   &  =T_{A}U_{\xi}\left(
T_{\operatorname*{Id}}R_{A}\left(  X\right)  \right)  =T_{\operatorname*{Id}%
}(U_{\xi}\circ R_{A})\left(  X\right)  =\\
&  =\left.  \frac{d}{dt}\right\vert _{t=0}f(\xi)f(\exp(Xt)A\xi)=f(\xi)\nabla
f^{\intercal}|_{A\xi}XA\xi,
\end{align*}
where $\nabla f^{\intercal}$ is the row vector $\left(  \partial f/\partial
\xi^{1},\ldots,\partial f/\partial\xi^{N}\right)  $. Now, let $\left\{
E_{ij}\right\}  _{0<i<j\leq N}$ be the orthonormal basis of $\mathfrak{so}(N)$
above-mentioned. By Section \ref{section Lie group methods} Equation
(\ref{antithetics 11}), the gradient $\nabla_{Y}(U_{\xi}\circ\varphi)$ in the
first part of (\ref{RMGradient}) is built from the composition of
$T_{\varphi(Y)}U_{\xi}^{sc}$ with $T_{Y}\varphi^{sc}$. Explicitly,
\[
\nabla_{Y}(U_{\xi}\circ\varphi)=\sum_{i<j}\left(  f(\xi)\nabla f^{\intercal
}|_{A\xi}T_{Y}\varphi^{sc}(E_{ij})A\xi\right)  E_{ij}\in\mathfrak{so}%
(N),~~A=\exp(Y)A_{0}.
\]
On the other hand, the second part of (\ref{RMGradient}) is simply minus the
gradient of the penalty function $P(Y)$, $-2Y\cdot\mathbf{1}_{B(0,\pi\sqrt
{R})^{c}}(Y)$. Thus, the Robbins-Monro part in the simulated annealing
algorithm is
\[
-\left(  \sum_{i<j}\left(  f(\xi)\nabla f^{\intercal}|_{A\xi}T_{Y}\varphi
^{sc}(E_{ij})A\xi\right)  E_{ij}\right)  -2Y\cdot\mathbf{1}_{B(0,\pi\sqrt
{R})^{c}}\left(  Y\right)
\]
evaluated at $Y=Y_{n-1}$, $A=A_{n-1}$, and $\xi=\xi_{n}$. This yields the
iterative scheme described in the statement.

To prove convergence we have to show that hypothesis \textbf{B1},\textbf{ B2
}and\textbf{ B3} in Theorem \ref{Theorem Fang-Gong-Quian} hold. Since
$U(\xi,A)$ is bounded the only part of the \textit{objective} function that
will contribute in the limits in those conditions will be the penalty function
$P(Y)$. That is, we can replace $\overline{F}(Y)$ by $-\nabla\Vert Y\Vert
^{2}=-2Y$. \textbf{B1},\textbf{ B2 }and\textbf{ B3} follow easily in this case.
\end{proof}

\begin{remark}
\normalfont \ 

\begin{enumerate}
\item The choice of initial matrix, $A_{0}$, is rather arbitrary unless one
has additional information on the behaviour of the function $f$.

However a couple of observations can be made to aid an informed decision. The
first is that $O(N)$ has two connected components, $SO(N)$ and $SO(N)^{-}$,
defined as the preimages of $+1$ and $-1$ by the continuous map $\mathrm{det}%
:O(N)\rightarrow\{ +1,-1 \}$. The algorithm we will define will stay in the
connected component that contains $A_{0}$.

The maximum of the covariance function $\operatorname*{Cov}\left(  f\left(
\xi\right)  ,f\left(  A\xi\right)  \right)  $ is achieved at
$A=\operatorname*{Id}\in SO(N)$, where it values $\operatorname*{Var}[f\left(
\xi\right)  ]$. Assuming the covariance function is continuous, it will be
positive around the identity matrix $\operatorname*{Id}$. Therefore, it seems
advisable to choose the initial matrix in the other connected component
$SO(N)^{-}$. It might be of interest to explore whether the minimum of the
covariance function always achieved in $SO(N)^{-}$.

\item The convergence of algorithm (\ref{algorithm steps}) depends on the
constants $d>0$. The parameter $d$ controls the \textit{heat injection} in the
annealing scheme. As mentioned in Theorem \ref{Theorem Fang-Gong-Quian}, $d$
must be larger than $d^{\ast}$ in Equation (\ref{antithetics 18}), which
measures the oscillatory nature of the function. In our case, given that the
objective function $\overline{U}$ is essentially the covariance, it can be
seen that $d^{\ast}\leq4\operatorname*{Var}[f\left(  \xi\right)  ]$. Since
$\operatorname*{Var}[f\left(  \xi\right)  ]$ is a number that any user of
Monte Carlo will be familiar with, we can estimate the magnitude of a valid
choice of $d$.

\item As stated the algorithm applies to functions $f$ which are
differentiable. Many of the functions used in pricing derivatives are non
differentiable along a hypersurface. In fact the algorithm will work provided
the set of non-differentiable points has measure zero a condition that will
hold in most real life scenarios.
\end{enumerate}
\end{remark}

\section{Examples\label{section examples}}

In this section, we are going to test the performance of our algorithm
(\ref{algorithm steps}) in finding the optimal antithetic for a couple of
elementary derivatives, an (arithmetic) Asian option and a covariance swap.
Both products will be priced under the Black-Scholes model. This means that,
under the risk-neutral probability, the price of an asset $S_{t}$ is expressed
as%
\begin{equation}
S_{t}=S_{0}\operatorname*{e}\nolimits^{(r-d-\nu^{2}/2)t+\nu B_{t}},
\label{antithetics 23}%
\end{equation}
where $r$ stands for the continuous interest rate, $d$ for the continuous
dividend yield, and $\nu$ for the volatility, all assumed constant; $B_{t}$
denotes a standard Brownian motion, i.e., $B_{t}-B_{s}\sim\mathcal{N}(0,t-s)$,
$t>s$. On the other hand, if $f_{t}(S_{1},...,S_{k})$ is the payoff function
at time $t$ of a given contract depending on $k$ assets, the price of such
contract at time $t$ is given by%
\begin{equation}
\operatorname*{e}\nolimits^{-rt}\operatorname*{E}\left[  f_{t}(S_{1}%
,...,S_{k})\right]  . \label{antithetics 24}%
\end{equation}

We do not intend to find the right price of such products in a realistic
market, where more sophisticated models may be required, but simply implement
algorithm (\ref{algorithm steps}) and show its efficiency and usefulness in
standard Monte Carlo simulations. We hope to explore antithetics under other
asset models in a future work.

\subsection{Asian option}

An {\bfseries\itshape(arithmetic) Asian call option} with strike price $K$ is
a derivative contract based on an asset price $S_{t}$ whose payoff at the
expiry time $T$ is given by
\[
\max\left(  0,\frac{1}{N}\sum_{i=1}^{N}S_{t_{i}}-K\right)  ,
\]
where $0\leq t_{1}\leq...\leq t_{N}=T$ is a sequence of times set in the
contract at which $S_{t}$ is observed. Consequently, according to
(\ref{antithetics 24}), the price of an Asian call option is%
\begin{equation}
\operatorname*{e}\nolimits^{-rT}\operatorname*{E}\left[  \max\left(
0,\frac{1}{N}\sum_{i=1}^{N}S_{t_{i}}-K\right)  \right]  .
\label{antithetics 25}%
\end{equation}
In order to estimate this expectation through a Monte Carlo simulation, we
replace each $S_{t_{i}}$ with the random variable%
\[
S_{t_{i}}=S_{t_{i-1}}\operatorname*{e}\nolimits^{(r-\nu^{2}/2)(t_{i}%
-t_{i-1})+\nu\sqrt{(t_{i}-t_{i-1})}\zeta},~\text{~}\zeta\sim\mathcal{N}(0,1),
\]
so that we can compute the sequence of prices iteratively. Therefore,
according to our picture, the expectation in (\ref{antithetics 25}) can be
rewritten as $\operatorname*{E}\left[  f(\xi)\right]  ,$ where%
\[
f(x_{1},...,x_{N})=\max\left(  0,\frac{1}{N}\sum_{i=1}^{N}%
%TCIMACRO{\dprod \limits_{j=1}^{i}}%
%BeginExpansion
{\displaystyle\prod\limits_{j=1}^{i}}
%EndExpansion
\operatorname*{e}\nolimits^{(r-\nu^{2}/2)(t_{j}-t_{j-1})+\nu\sqrt
{(t_{j}-t_{j-1})}x_{j}}-K\right)  \text{ and }\xi\sim\mathcal{N}%
(0,\operatorname*{Id}\nolimits_{N}).
\]
In the previous expression, we implicitly assumed $t_{0}=0$.

In this particular example, we are going to find the optimal antithetic for an
Asian option averaged over the asset prices on the first of each month along
one year. For the sake of simplicity, we suppose that all months have the same
number of days. That is, the expiry time $T=1$ equals one year and
$t_{i}=i/12$ years, $i=1,...,12$. More explicitly, our antithetic matrix will
be of dimension $12$, which means that our optimizing problem
(\ref{antithetics 4}) takes place on a space of dimension $\dim(\mathfrak{so}%
(12))=66$. We will suppose that the asset price at the initial date is
$S_{0}=\pounds 100$ and that the strike price equals $\pounds 100$ as well.
Furthermore, we set $r=2.83\%$ and $\nu=10.36\%$.

The results obtained are summarized in Table 1. The first row refers to a
crude Monte Carlo estimation of one Asian option according to the data
above-mentioned. The second and the third rows contain, respectively,
information on the Monte Carlo estimation of (\ref{antithetics 24}) using the
antithetics procedure (\ref{antithetics 26}) when the matrix $A$ is minus the
identity $-\operatorname*{Id}_{12}$ and the antithetic $A^{\ast}$ provided by
algorithm (\ref{algorithm steps}) with $\gamma=1/2$ and initial condition
$A_{0}=-\operatorname*{Id}_{12}$. For any estimation, we specify the price of
the option, the estimated variance $\sigma^{2}=\operatorname*{Var}%
[\operatorname*{e}\nolimits^{-rT}f\left(  \xi\right)  ]$ of the method, the
error $\sigma/\sqrt{n}$ associated to that estimation, and the time (in
seconds) our computer (a laptop with a 1.6 MHz CPU) invested in those
computations.%
\[
\fbox{$%
\begin{array}
[c]{rcccc}
& \text{Price (\pounds )} & \text{Variance }\sigma^{2} & \text{MC Error
(}\sigma/\sqrt{n}\text{)} & \text{Time (s)}\\
\text{Crude Monte Carlo} & 3.15 & 17.25 & 0.01 & 8.75\\
\text{Antithetic }-\operatorname*{Id}_{12} & 3.15 & 3.50 & 0.01 & 3.41\\
\text{Antithetic }A^{\ast} & 3.15 & 3.60 & 0.01 & 3.54
\end{array}
$}%
\]

\begin{center}
{\small Table 1: Monte Carlo simulations for an Asian option\medskip}
\end{center}

First of all, we observe that the use of antithetics dramatically reduces the
variance a crude Monte Carlo simulation by, roughly speaking, $5$, which
means, once an error threshold is fixed, we need less than half of the time
required in a crude Monte Carlo to price the option with the same accuracy.
However, we can see that the use of antithetics in the original Hammersley
sense, that is, taking $-\operatorname*{Id}_{12}$, slightly beats the
antithetic provided by algorithm (\ref{algorithm steps}). This is because, in
this particular example, $-\operatorname*{Id}_{12}$ is one of the points where
the covariance function $\operatorname*{Cov}\left(  f(\xi),f(A\xi)\right)  $
reaches its minimum value. If algorithm (\ref{algorithm steps}) does not
provide an antithetic so performing is because we can only run it for a finite
time while its convergence to the set $\underline{S}=\min_{A\in\mathfrak{so}%
(12)}\operatorname*{Cov}\left(  f(\xi),f(A\xi)\right)  $ is assured as
$t\rightarrow\infty$. However, despite this fact and that $A^{\ast}$ is
obtained by a numeric method, which by definition intrinsically carries some
degree of approximation, the result is satisfactory and we conclude that
$A^{\ast}$ is very close to an absolute minimum.

The main drawback of algorithm (\ref{algorithm steps}) is its speed of
convergence. Since at each step we have to compute, on the one hand, the
exponential of a skew-symmetric matrix (of dimension $12$) and, on the other,
the tangent map $T\exp^{sc}$ to write down the gradient
(\ref{algorithm step 1}), the algorithm is very time-consuming. For example,
the antithetic $A^{\ast}$ used in the simulations of Table 1 is obtained after
$10000$ iterations, which our 1.6MHz computer carried out in approximately
forty minutes, nothing comparable with the times spent in the Monte Carlo
simulations (see Table 1). Therefore, finding the optimal antithetic $A^{\ast
}$ before a crude Monte Carlo seems advisable if $A^{\ast}$ may be reused in
several simulations. We hope to improve the efficiency of the algorithm in a
future work.

\subsection{Covariance swap}

The next example aims at showing that the optimal antithetic $A^{\ast}$ of
some daily traded products can be very different from the antithetic
$-\operatorname*{Id}$ which, moreover, turns out to be completely useless.

A {\bfseries\itshape covariance swap} depending on two assets $S_{1}$ and
$S_{2}$ is a contract that, at the expiry date $T$, pays%
\begin{equation}
\frac{1}{N}\sum_{i=1}^{N}\left(  \frac{(S_{1})_{t_{i}}}{(S_{1})_{t_{i-1}}%
}-1\right)  \left(  \frac{(S_{2})_{t_{i}}}{(S_{2})_{t_{i-1}}}-1\right)
\label{antithetics 27}%
\end{equation}
where, as in the case of the Asian option, $0=t_{0}\leq t_{1}\leq...\leq
t_{N}=T$ is a sequence of times fixed in the contract. This quantity can be
negative, which means that the holder of the contract, instead of receiving
any money, must pay that quantity. It is called covariance swap because,
roughly speaking, it measures the realised covariance between the returns of
the assets $S_{1}$ and $S_{2}$. Sometimes, the returns are alternatively
expressed as $\ln\left(  S_{t_{i}}/S_{t_{i-1}}\right)  $.

In the Balck-Scholes model, where the assets are uncorrelated, the payoff
(\ref{antithetics 27}) is very small (but not zero though, because the spot
prices are only observed on a finite number of dates and not
\textit{continuously}). In order to increase the price of
(\ref{antithetics 27}) so that we avoid working with too many decimal places,
we take the maximum between (\ref{antithetics 27}) and $0$,%
\begin{equation}
\max\left(  0,\frac{1}{N}\sum_{i=1}^{N}\left(  \frac{(S_{1})_{t_{i}}}%
{(S_{1})_{t_{i-1}}}-1\right)  \left(  \frac{(S_{2})_{t_{i}}}{(S_{2})_{t_{i-1}%
}}-1\right)  \right)  . \label{antithetics 28}%
\end{equation}
We insist, however, that this is not the way covariance swaps are traded in
real markets, where they are mainly used as a hedge against possible unwanted
correlation effects between assets or indices.

We are going to fix the details of our simulation. Let $\nu_{1}=17.36\%$ and
$\nu_{2}=10.12\%$ be the volatilities of $S_{1}$ and $S_{2}$. Suppose that the
asset $S_{2}$ pays a constant dividend yield $d=2.03\%$ and that the constant
free-risk interest rate is again $r=2.83\%$. Hence the assets $S_{1}$ and
$S_{2}$ follow the log-normal processes:%

\[
(S_{1})_{t}=\left(  S_{1}\right)  _{0}\operatorname*{e}\nolimits^{(r-\nu
_{1}^{2}/2)t+\nu_{1}B_{t}^{(1)}},~~(S_{2})_{t}=\left(  S_{2}\right)
_{0}\operatorname*{e}\nolimits^{(r-d-\nu_{2}^{2}/2)t+\nu_{2}B_{t}^{(2)}},
\]
where $B_{t}^{(1)}$ and $B_{t}^{(2)}$ are two independent Brownian motions. As
in the previous example, the price of (\ref{antithetics 28}) can be written as
$\operatorname*{e}^{-rT}\operatorname*{E}\left[  f(\xi)\right]  $. Explicitly,%
\begin{gather}
f(x_{1},...,x_{2N})=\max\left(  0,\frac{1}{N}\sum\nolimits_{j=1}^{N}\left(
\operatorname*{e}\nolimits^{(r-\nu_{1}^{2}/2)(t_{j}-t_{j-1})+\nu_{1}%
\sqrt{(t_{j}-t_{j-1})}x_{2j-1}}-1\right)  \cdot\right. \nonumber\\
\hspace{5.3cm}\left.  \cdot\left(  \operatorname*{e}\nolimits^{(r-d-\nu
_{2}^{2}/2)(t_{j}-t_{j-1})+\nu_{2}\sqrt{(t_{j}-t_{j-1})}x_{2j}}-1\right)
\right)  \label{antithetics 29}%
\end{gather}
and $\xi\sim\mathcal{N}(0,\operatorname*{Id}\nolimits_{2N})$ is a Gaussian
vector whose odd components $\xi^{2j+1}$, $j=1,...,N$, are used to generate
the evolution of $S_{1}$ while the even components $\xi^{2j}$ go with $S_{2}$.
It is clear from (\ref{antithetics 29}) that the values of $S_{1}$ and $S_{2}$
at time $t=0$ play no role in the final price of the covariance swap option,
so we can take them equal to $1$.

Table 2 summarises the prices of $100$ covariance swaps estimated by Monte
Carlo. As for the Asian option before studied, the asset prices $S_{t_{i}}$
are read on the first of every month along one year, where all months are
supposed to have $30$ days, which means that the expiry time $T=1$ equals one
year and $t_{i}=i/12$ years, $i=1,...,12$. Our (optimal) antithetic matrix
will be of dimension $24$, so that our optimizing problem (\ref{antithetics 4}%
) takes place on a space of dimension $\dim(\mathfrak{so}(24))=276$. As in
Table 1, the first row refers to a crude Monte Carlo, the second one to the
antithetics procedure (\ref{antithetics 26}) when the matrix $A$ is minus the
identity $-\operatorname*{Id}_{24}$, and the third one to the antithetic
$A^{\ast}$ provided by algorithm (\ref{algorithm steps}) with $\gamma=1/2$ and
initial condition $A_{0}=-\operatorname*{Id}_{24}$. Again, we specify the
price of the option, the estimated variance $\sigma^{2}=\operatorname*{Var}%
[\operatorname*{e}\nolimits^{-rT}f\left(  \xi\right)  ]$ of the method, the
error $\sigma/\sqrt{n}$ associated to that estimation, and the time (in
seconds) spent in those simulations.%
\[
\fbox{$%
\begin{array}
[c]{rcccc}
& \text{Price (\pounds )} & \text{Variance }\sigma^{2} & \text{MC Error
(}\sigma/\sqrt{n}\text{)} & \text{Time (s)}\\
\text{Crude Monte Carlo} & 0.198 & 0.081 & 0.001 & 7.28\\
\text{Antithetic }-\operatorname*{Id}_{24} & 0.197 & 0.081 & 0.001 & 15.72\\
\text{Antithetic }A^{\ast} & 0.198 & 0.029 & 0.001 & 5.02
\end{array}
$}%
\]

\begin{center}
{\small Table 2: Monte Carlo simulations for an option on a covariance
swap\medskip}
\end{center}

We can observe that now the antithetic $A=-\operatorname*{Id}_{24}$ does not
improve the efficiency of the crude Monte Carlo. Indeed, it is easy to check
that, in this particular example, $\operatorname*{Cov}(f(\xi),f(-\xi
))=\operatorname*{Var}[f(\xi)]$, $\xi\sim\mathcal{N}(0,\operatorname*{Id}%
\nolimits_{24})$, so that the variance of (\ref{antithetics 26}) is the same
of the original (crude) Monte Carlo. However, since at each step we now
evaluate the payoff function $f$ twice, we need twice as much time. In
conclusion, the antithetic $A=-\operatorname*{Id}_{24}$ is useless.

On the contrary, the antithetic $A^{\ast}$ reduces the variance $2.8$ times
which, in turn, means that the computational time is reduced, approximately,
by $1.45$. Unfortunately, unlike $-\operatorname*{Id}_{24}$, $A^{\ast}$ is an
orthogonal dense matrix that has no trivial interpretation. And what is worse,
as the dimensionality of the optimisation problem grows, so does the time
needed to obtain a good approximation of $A^{\ast}$. For example, in the
location of the skew-symmetric matrix $Y^{\ast}\in\mathfrak{so}(24)$ such that
$\exp\left(  Y^{\ast}\right)  =A^{\ast}$, a 1.6MHz laptop invested more than
an hour in carrying out $10000$ iterations of (\ref{algorithm steps}).
Although there are certainly more performing computers, this time is huge
compared with the times of the simulations we want to improve (see Table 2).

\section{A dynamically optimized Monte Carlo\label{section central limit}}

Section \ref{section antithetics revisited} provides an algorithm to locate
the optimal antithetic $A^{\ast}$. Recall that the use of this matrix
$A^{\ast}$ will improve the calculation of $m=\operatorname*{E}[f(\xi)]$ as a
Monte Carlo average in the form
\[
\lim_{n\rightarrow\infty}\frac{1}{n}\sum_{k=1}^{n}\frac{f(\xi_{k})+f(A^{\ast
}\xi_{k})}{2}.
\]
In practice, as we saw in Section \ref{section examples}, it seems too
expensive to locate $A^{\ast}$ and then use it in the Monte Carlo calculation.
It would be better to use the antithetics dynamically, meaning that we use as
approximating sum
\begin{equation}
\lim_{n\rightarrow\infty}\frac{1}{n}\sum_{k=1}^{n}\frac{f(\xi_{k}%
)+f(A_{k-1}\xi_{k})}{2}.\label{antithetics 19}%
\end{equation}
In other words, we calculate the optimal antithetic at the same time as we
estimate $\operatorname*{E}[f(\xi)]$.

It is worth emphasizing that the random samples $\{\zeta_{n}\}_{n\in
\mathbb{N}}$\ used in the location of $A^{\ast}$ via the simulated annealing
(\ref{algorithm steps}) can be also recycled in (\ref{antithetics 19}) to
approximate $m$. Did we do that, the results in this section would continue
being valid. However, the expressions would become more complex but with no
additional content. Therefore, in order to illustrate that antithetic variates
can be used dynamically, we are only going to study the estimations of $m$
obtained from (\ref{antithetics 19}). Furthermore, for the benefit of a
simpler notation, we will abbreviate%
\[
g(A_{k-1},\xi_{k}):=\frac{1}{2}\left(  f(\xi_{k})+f(A_{k-1}\xi_{k})\right)  .
\]

The random variables $\left\{  g\left(  A_{k-1},\xi_{k}\right)  \right\}
_{k\in\mathbb{N}}$ are neither independent nor equally distributed because,
for any $k\in\mathbb{N}$, $A_{k}$ is not deterministic but depends on the
previous random variables $\left\{  \xi_{i},\zeta_{i}\right\}  _{i=1,...,k}$.
Consequently, we cannot invoke the Strong Law of Large Numbers to guarantee
that (\ref{antithetics 19}) converges to $m$ almost surely. In other words,
(\ref{antithetics 19}) is different from computing $m$ as%
\[
\frac{1}{n}\sum_{k=1}^{n}\frac{f\left(  \xi_{k}\right)  +f\left(  A^{\ast}%
\xi_{k}\right)  }{2}%
\]
for the optimal value $A^{\ast}$, where $A^{\ast}$ is supposed to be
deterministic and fixed. Moreover, the dependence of $\left\{  g\left(
A_{k-1},\xi_{k}\right)  \right\}  _{k\in\mathbb{N}}$ might cause the
appearance of some positive correlations which might spoil the efficiency of
our method. Nevertheless, our situation is not that bad. Indeed, if we define
\[
S_{n}:=\frac{1}{n}\sum_{k=1}^{n}g\left(  A_{k-1},\xi_{k}\right)  ,
\]
we are going to prove in Theorem \ref{prop 1} that, on the one hand, $S_{n}$
converges to the expectation $m=\operatorname*{E}[f\left(  \xi\right)  ]$
almost surely as $n\rightarrow\infty$ and, on the other, that $\sqrt{n}%
(S_{n}-m)$ converges in law to a normal variable $\mathcal{N}\left(
0,\sigma_{\ast}^{2}\right)  $ with variance%
\[
\sigma_{\ast}^{2}=\operatorname*{Var}\left[  g\left(  A^{\ast},\xi\right)
\right]  =\frac{1}{2}\left(  \operatorname*{Var}[f(\xi)]+\operatorname*{Cov}%
\left(  f\left(  \xi\right)  ,f\left(  A^{\ast}\xi\right)  \right)  \right)
,
\]
$\xi\sim\mathcal{N}(0,\operatorname*{Id}\nolimits_{N})$. It is worth noticing
that $\sigma_{\ast}^{2}=\operatorname*{Var}\left[  g\left(  A^{\ast}%
,\xi\right)  \right]  $ is the smallest possible variance we can get using
antithetic variates.

Finally, in order to estimate $\sigma_{\ast}^{2}$ and thus obtain empirical
confidence intervals for $m$, we prove in Theorem \ref{central limit theorem}
that%
\[
\frac{1}{n}\sum_{k=1}^{n}g^{2}\left(  A_{k-1},\xi_{k}\right)  -S_{n}%
^{2}\underset{n\rightarrow\infty}{\longrightarrow}\sigma_{\ast}^{2}\text{ \ in
probability.}%
\]
The following results are inspired by \cite{arouna 1, arouna 2}.

\begin{theorem}
\label{prop 1}Let $A^{\ast}\in O(N)$ and suppose that there exists a sequence
$\left\{  A_{k}\right\}  _{k\in\mathbb{N}}$ of random variables taking values
in $O(N)$ such that $A_{k}\rightarrow A^{\ast}$ in probability as
$k\rightarrow\infty$. Let $\left\{  \xi_{k}\right\}  _{k\in\mathbb{N}}$ and
$\left\{  \zeta_{k}\right\}  _{k\in\mathbb{N}}$ be a couple of sequences of
i.i.d Gaussian vectors mutually independent. Let $\mathcal{F}_{k}%
=\sigma\left(  \xi_{i},\zeta_{i}~|~1\leq i\leq k\right)  $ be the $\sigma
$-algebra generated by $\{\xi_{i},\zeta_{i}\}_{i=1,...,k}$ and suppose that
$A_{k}$ is $\mathcal{F}_{k}$-measurable. Furthermore, if $\xi\sim
\mathcal{N}(0,\operatorname*{Id}\nolimits_{N})$ denote an independent Gaussian
vector, assume that $\operatorname*{E}[|g\left(  A_{k-1},\xi\right)
|^{2}]<\infty$ for any $k\in\mathbb{N}$ and that the map $A\mapsto
\operatorname*{E}[\left\vert g\left(  A,\xi\right)  \right\vert ^{2}]$ is
continuous at $A^{\ast}$. Then,%
\[
\frac{1}{n}\sum_{k=1}^{n}g\left(  A_{k-1},\xi_{k}\right)  \underset
{n\rightarrow\infty}{\overset{a.s.}{\longrightarrow}}m.
\]

\end{theorem}

\begin{proof}
Define $M_{0}=0$ and%
\[
M_{n}=\sum_{k=1}^{n}\left(  g\left(  A_{k-1},\xi_{k}\right)  -m\right)
,~n\geq1\text{.}%
\]
Let $\mathcal{F}_{n}$ be the $\sigma$-algebra generated by $\left\{  \xi
_{i},\zeta_{i}\right\}  _{i=1,...,n}$. Since $A_{n-1}$ is $\mathcal{F}_{n-1}%
$-measurable and $\xi_{n}$ is independent of $\mathcal{F}_{n-1}$, we have%
\[
\operatorname*{E}\left[  g\left(  A_{n-1},\xi_{n}\right)  ~|~\mathcal{F}%
_{n-1}\right]  =\left.  \operatorname*{E}\left[  g\left(  A,\xi_{n}\right)
\right]  \right\vert _{A=A_{n-1}}=m
\]
and $\left\{  M_{n}\right\}  _{n\in\mathbb{N}}$ is a martingale with respect
to the filtration $\left\{  \mathcal{F}_{n}\right\}  _{n\in\mathbb{N}}$. As we
imposed that $\operatorname*{E}[\left(  g\left(  A_{n-1},\xi_{n}\right)
\right)  ^{2}]<\infty$ for any $n\in\mathbb{N}$, it is not difficult to see
that $\left\{  M_{n}\right\}  _{n\in\mathbb{N}}$ is a square integrable
martingale. On the other hand, its angle bracket is given by%
\begin{align*}
\left\langle M\right\rangle _{n}  &  =\sum_{k=1}^{n}\operatorname*{E}\left[
\left(  \Delta M_{k}\right)  ^{2}~|~\mathcal{F}_{k-1}\right]  =\sum_{k=1}%
^{n}\left(  \operatorname*{E}\left[  g^{2}\left(  A_{k-1},\xi_{k}\right)
~|~\mathcal{F}_{k-1}\right]  -m^{2}\right) \\
&  =\sum_{k=1}^{n}\left(  \left.  \operatorname*{E}\left[  g^{2}\left(
A,\xi_{k}\right)  \right]  \right\vert _{A=A_{k-1}}-m^{2}\right)  .
\end{align*}
Now, since $A_{k}\rightarrow A^{\ast}$ in probability as $k\rightarrow\infty$
and $A\mapsto\operatorname*{E}[g^{2}(A,\xi_{k})]$ is continuous at $A^{\ast}$,
$\left.  \operatorname*{E}[g^{2}(A,\xi_{k})]\right\vert _{A=A_{k-1}}$ also
converges in probability to $\left.  \operatorname*{E}[g^{2}(A,\xi
_{k})]\right\vert _{A=A^{\ast}}$ by Lemma (\ref{lema continuitat}). Applying
Lemma (\ref{Cesaro Lemma}), we see that%
\[
\frac{\left\langle M\right\rangle _{n}}{n}\underset{n\rightarrow\infty
}{\longrightarrow}\left.  \operatorname*{E}\left[  g^{2}\left(  A,\xi\right)
\right]  \right\vert _{A=A^{\ast}}-m^{2}=\operatorname*{Var}\left[  g\left(
A^{\ast},\xi\right)  \right]
\]
in probability. Therefore, by Theorem \ref{teorema auxiliar} (i), we conclude
that $\frac{1}{n}M_{n}\underset{n\rightarrow\infty}{\longrightarrow}0$ a.s.
\end{proof}

\begin{theorem}
[Central Limit Theorem]\label{central limit theorem} With the same notation as
in Theorem \ref{prop 1}, assume now that $\operatorname*{E}[|g\left(
A_{k-1},\xi\right)  |^{4}]<\infty$ for any $k\in\mathbb{N}$ and that
$A\mapsto\operatorname*{E}[\left\vert g\left(  A,\xi\right)  \right\vert
^{p}]$ is continuous at $A^{\ast}$, $1\leq p\leq4$.

\begin{description}
\item[(i)] The following convergence in law holds:%
\[
\sqrt{n}\left(  S_{n}-m\right)  \underset{n\rightarrow\infty}{\longrightarrow
}\mathcal{N}\left(  0,\sigma_{\ast}^{2}\right)  ,
\]
where $\sigma_{\ast}^{2}:=\operatorname*{Var}\left[  g\left(  A^{\ast}%
,\xi\right)  \right]  $.

\item[(ii)] Let $\sigma_{n}^{2}:=\frac{1}{n}\sum_{k=1}^{n}g^{2}\left(
A_{k-1},\xi_{k}\right)  -S_{n}^{2}$. Then $\sigma_{n}^{2}\overset
{n\rightarrow\infty}{\longrightarrow}\sigma_{\ast}^{2}$ in probability.
\end{description}
\end{theorem}

\begin{proof}
\ \ \ 

\begin{description}
\item[(i)] First of all, observe that%
\begin{align*}
\operatorname*{E}\left[  |g\left(  A_{k-1},\xi_{k}\right)  -m|^{4}%
~|~\mathcal{F}_{k-1}\right]   &  =\operatorname*{E}\left[  g^{4}\left(
A_{k-1},\xi_{k}\right)  ~|~\mathcal{F}_{k-1}\right]  -3m^{4}\\
&  -4m\operatorname*{E}\left[  g^{3}\left(  A_{k-1},\xi_{k}\right)
~|~\mathcal{F}_{k-1}\right]  +6m^{2}\operatorname*{E}\left[  g^{2}\left(
A_{k-1},\xi_{k}\right)  ~|~\mathcal{F}_{k-1}\right] \\
&  =\left.  \operatorname*{E}\left[  g^{4}\left(  A,\xi_{k}\right)  \right]
\right\vert _{A=A_{k-1}}-3m^{4}\\
&  -4m\left.  \operatorname*{E}\left[  g^{3}\left(  A,\xi_{k}\right)  \right]
\right\vert _{A=A_{k-1}}+6m^{2}\left.  \operatorname*{E}\left[  g^{2}\left(
A,\xi_{k}\right)  \right]  \right\vert _{A=A_{k-1}}.
\end{align*}
Since $A\mapsto\operatorname*{E}[|g\left(  A_{k-1},\xi\right)  |^{p}]$ is
continuous and $A_{k}\rightarrow A$ in probability as $k\rightarrow\infty$, by
Lemma (\ref{lema continuitat}) we have that $\operatorname*{E}\left[
|g\left(  A_{k-1},\xi_{k}\right)  -m|^{4}~|~\mathcal{F}_{k-1}\right]  $
converges in probability to a positive random variable $L\in L_{\mathbb{R}%
}^{0}\left(  \Omega,P\right)  $ as $k\rightarrow\infty$. Thus, by Lemma
(\ref{Cesaro Lemma}),%
\[
\frac{1}{n}\sum_{k=1}^{n}\operatorname*{E}\left[  |g\left(  A_{k-1},\xi
_{k}\right)  -m|^{4}~|~\mathcal{F}_{k-1}\right]  \underset{n\rightarrow\infty
}{\longrightarrow}L
\]
in probability.

For any $a>0$, define%
\[
F_{n}:=\frac{1}{n}\sum_{k=1}^{n}\operatorname*{E}\left[  |g\left(  A_{k-1}%
,\xi_{k}\right)  -m|^{2}\mathbf{1}_{\{|g\left(  A_{k-1},\xi_{k}\right)
-m|>a\}}~|~\mathcal{F}_{k-1}\right]  .
\]
It can be easily checked that%
\[
F_{n}\left(  a\right)  \leq\frac{1}{na^{2}}\sum_{k=1}^{n}\operatorname*{E}%
\left[  |g\left(  A_{k-1},\xi_{k}\right)  -m|^{4}~|~\mathcal{F}_{k-1}\right]
\]
so that, taking the limit superior of the sequence $\left\{  F_{n}\left(
a\right)  \right\}  _{n\in\mathbb{N}}$ in probability, we obtain that
$\lim\sup_{n\rightarrow\infty}F_{n}\left(  a\right)  \leq a^{-2}L$. If now
$a=\varepsilon\sqrt{n}$ with $\varepsilon>0$ fixed, we have
\[
\underset{n\rightarrow\infty}{\lim\sup}F_{n}\left(  \varepsilon\sqrt
{n}\right)  =0
\]
in probability and the Lindberg's condition holds. Therefore, $\sqrt{n}\left(
S_{n}-m\right)  \overset{n\rightarrow\infty}{\longrightarrow}\mathcal{N}%
\left(  0,\sigma_{\ast}^{2}\right)  $ in law by Theorem \ref{teorema auxiliar} (ii).

\item[(ii)] Essentially, we only need to show that $\frac{1}{n}\sum_{k=1}%
^{n}\left.  \operatorname*{E}\left[  g^{2}\left(  A,\xi\right)  \right]
\right\vert _{A=A_{k-1}}$ and $\frac{1}{n}\sum_{k=1}^{n}g^{2}\left(
A_{k-1},\xi_{k}\right)  $ have the same limit a.s.. We are going to do so by
mimicking the proof of Theorem \ref{prop 1} and quoting Theorem
\ref{teorema auxiliar}.

Define $M_{0}=0$ and%
\[
M_{n}:=\sum_{k=1}^{n}\left(  g^{2}\left(  A_{k-1},\xi_{k}\right)  -\left.
\operatorname*{E}\left[  g^{2}\left(  A,\xi_{k}\right)  \right]  \right\vert
_{A=A_{k-1}}\right)  ,~n\geq1.
\]
It is not difficult to check that $\left\{  M_{n}\right\}  _{n\in\mathbb{N}}$
is a square integrable martingale with respect to $\left\{  \mathcal{F}%
_{n}\right\}  _{n\in\mathbb{N}}$, $\mathcal{F}_{n}=\sigma\left(  \xi_{i}%
,\zeta_{i}~|~i=1,...,n\right)  $, whose angle bracket is given by%
\begin{align*}
\left\langle M\right\rangle _{n}  &  =\sum_{k=1}^{n}\left(  \operatorname*{E}%
\left[  g^{4}\left(  A_{k-1},\xi_{k}\right)  ~|~\mathcal{F}_{k-1}\right]
-\left.  \operatorname*{E}\left[  g^{2}\left(  A,\xi_{k}\right)  \right]
^{2}\right\vert _{A=A_{k-1}}\right) \\
&  =\sum_{k=1}^{n}\left(  \left.  \operatorname*{E}\left[  g^{4}\left(
A,\xi_{k}\right)  \right]  \right\vert _{A=A_{k-1}}-\left.  \operatorname*{E}%
\left[  g^{2}\left(  A,\xi_{k}\right)  \right]  ^{2}\right\vert _{A=A_{k-1}%
}\right)  .
\end{align*}
By the continuity of $A\mapsto\operatorname*{E}\left[  \left\vert g\left(
A,\xi\right)  \right\vert ^{p}\right]  $, $1\leq p\leq4$, and Lemma
\ref{Cesaro Lemma} we have%
\[
\frac{1}{n}\left\langle M\right\rangle _{n}\underset{n\rightarrow\infty
}{\longrightarrow}\left.  \operatorname*{E}\left[  g^{4}\left(  A,\xi\right)
\right]  \right\vert _{A=A^{\ast}}-\left.  \operatorname*{E}\left[
g^{2}\left(  A,\xi\right)  \right]  ^{2}\right\vert _{A=A^{\ast}%
}=\operatorname*{Var}\left[  g^{2}\left(  A^{\ast},\xi\right)  \right]
\]
in probability. Then Theorem \ref{teorema auxiliar} (i) claims that%
\[
\frac{1}{n}\sum_{k=1}^{n}\left.  \operatorname*{E}\left[  g^{2}\left(
A,\xi_{k}\right)  \right]  \right\vert _{A=A_{k-1}}-\frac{1}{n}\sum_{k=1}%
^{n}g^{2}\left(  A_{k-1},\xi_{k}\right)  \underset{n\rightarrow\infty
}{\longrightarrow}0\text{ a.s..}%
\]
Since $\left.  \operatorname*{E}\left[  g^{2}\left(  A,\xi_{k}\right)
\right]  \right\vert _{A=A_{k-1}}=\left.  \operatorname*{E}\left[
g^{2}\left(  A,\xi\right)  \right]  \right\vert _{A=A_{k-1}}$ for any
arbitrary Gaussian vector $\xi$, we conclude that $\frac{1}{n}\sum_{k=1}%
^{n}g^{2}\left(  A_{k-1},\xi_{k}\right)  $ and $\frac{1}{n}\sum_{k=1}%
^{n}\left.  \operatorname*{E}\left[  g^{2}\left(  A,\xi\right)  \right]
\right\vert _{A=A_{k-1}}$ have the same limit almost surely.

Finally, as $\frac{1}{n}\sum_{k=1}^{n}\left.  \operatorname*{E}\left[
g^{2}\left(  A,\xi\right)  \right]  \right\vert _{A=A_{k-1}}$ and $S_{n}^{2}$
converge in probability to $\operatorname*{E}\left[  g^{2}\left(  A^{\ast}%
,\xi\right)  \right]  $ and $\operatorname*{E}\left[  g\left(  A^{\ast}%
,\xi\right)  \right]  ^{2}$ as $n\rightarrow\infty$ respectively, we obtain
that%
\[
\sigma_{n}^{2}=\frac{1}{n}\sum_{k=1}^{n}g^{2}\left(  A_{k-1},\xi_{k}\right)
-S_{n}^{2}\underset{n\rightarrow\infty}{\longrightarrow}\operatorname*{E}%
\left[  g^{2}\left(  A^{\ast},\xi\right)  \right]  -\operatorname*{E}\left[
g\left(  A^{\ast},\xi\right)  \right]  ^{2}=\operatorname*{Var}\left[
g\left(  A^{\ast},\xi\right)  \right]
\]
in probability.
\end{description}
\end{proof}

\appendix

\section{Appendix}

In this appendix we recall some auxiliary results. The first one is a
probabilistic version of the so-called Cesaro Lemma and reads

\begin{lemma}
\label{Cesaro Lemma}Let $\left\{  \xi_{n}\right\}  _{n\in\mathbb{N}}\subset
L_{\mathbb{R}^{r}}^{0}\left(  \Omega,P\right)  $ be a sequence of random
variables variables and let $\xi\in L_{\mathbb{R}^{r}}^{0}\left(
\Omega,P\right)  $. If $\xi_{n}\rightarrow\xi$ in probability as
$n\rightarrow\infty$, then%
\[
\frac{1}{n}\sum_{k=1}^{n}\xi_{k}\underset{n\rightarrow\infty}{\longrightarrow
}\xi\text{ in probability as well.}%
\]

\end{lemma}

\begin{proof}
Let $\zeta$, $\eta\in L_{\mathbb{R}^{r}}^{0}\left(  \Omega,P\right)  $. The
function $d\left(  \zeta,\eta\right)  =\operatorname*{E}\left[  \min\left(
1,\left\Vert \zeta-\eta\right\Vert \right)  \right]  $ is a distance in
$L_{\mathbb{R}^{r}}^{0}\left(  \Omega,P\right)  $ compatible with convergence
in probability. That is, a sequence $\left\{  \eta_{n}\right\}  _{n\in
\mathbb{N}}\subset L_{\mathbb{R}^{r}}^{0}\left(  \Omega,P\right)  $ converges
in probability to $\eta\in L_{\mathbb{R}^{r}}^{0}\left(  \Omega,P\right)  $ if
and only if $d\left(  \eta_{n},\eta\right)  \rightarrow0$ as $n\rightarrow
\infty$.

Let $\varepsilon>0$ and let $n_{0}\in\mathbb{N}$ be large enough so that
$d(\xi_{n},\xi)\leq\varepsilon/2$ for any $n>n_{0}$. We can choose
$n_{1}>n_{0}$ so that%
\[
\frac{1}{n_{1}}\sum_{k=1}^{n_{0}}d\left(  \xi_{k},\xi\right)  \leq
\frac{\varepsilon}{2}.
\]
Then, for any $n>n_{1}$ we have
\begin{align*}
d\left(  \frac{1}{n}\sum\nolimits_{k=1}^{n}\xi_{k},\xi\right)   &  \leq
\frac{1}{n}\sum\nolimits_{k=1}^{n}d\left(  \xi_{k},\xi\right)  =\frac{1}%
{n}\sum\nolimits_{k=1}^{n_{0}}d\left(  \xi_{k},\xi\right)  +\frac{1}{n}%
\sum\nolimits_{k=n_{0}}^{n_{1}}d\left(  \xi_{k},\xi\right) \\
&  \leq\frac{1}{n_{1}}\sum\nolimits_{k=1}^{n_{0}}d\left(  \xi_{k},\xi\right)
+\frac{\varepsilon(n-n_{0})}{2n}\leq\varepsilon.
\end{align*}
Consequently, $\frac{1}{n}\sum_{k=1}^{n}\xi_{k}$ converges to $\xi$ in
probability.\smallskip
\end{proof}

The next lemma is rather elementary and well known. However, we state and
prove it for the benefit of a clearer exposition.

\begin{lemma}
\label{lema continuitat}Let $X$ be a topological space. Let $\left\{
x_{n}\right\}  _{n\in\mathbb{N}}\subset L_{X}^{0}\left(  \Omega,P\right)  $ be
a sequence of random variables variables and $x\in X$. Let $f:X\rightarrow
\mathbb{R}$ be a function continuous at $x$. If $x_{n}\rightarrow x$ in
probability as $n\rightarrow\infty$, then $f(x_{n})\rightarrow f\left(
x\right)  $ in probability as well.
\end{lemma}

\begin{proof}
Let $\varepsilon>0$. By the continuity of $f$, there exists an open
neighborhood $V_{x}$ around $x$ such that $\left\vert f\left(  y\right)
-f(x)\right\vert \leq\varepsilon$ for any $y\in V_{x}$. Then,%
\[
\left\{  \omega\in\Omega~|~\left\vert f\left(  x_{n}\right)  -f(x)\right\vert
>\varepsilon\right\}  \subseteq\left\{  \omega\in\Omega~|~x_{n}\notin
V_{x}\right\}  .
\]
Since $P\left(  \left\{  \omega\in\Omega~|~x_{n}\notin V_{x}\right\}  \right)
\rightarrow0$ as $n\rightarrow\infty$ because $\left\{  x_{n}\right\}
_{n\in\mathbb{N}}$ converges in probability to $x\in\mathbb{R}$, we see that%
\[
P\left(  \left\{  \omega\in\Omega~|~\left\vert f\left(  x_{n}\right)
-f(x)\right\vert >\varepsilon\right\}  \right)  \underset{n\rightarrow\infty
}{\longrightarrow}0
\]
as well, and $\left\{  f(x_{n})\right\}  _{n\in\mathbb{N}}$ converges in
probability to $f\left(  x\right)  $.\smallskip
\end{proof}

Unlike the previous results, the following theorem is much deeper. Roughly
speaking, it is a powerful tool to prove general versions of the Law of Large
Numbers and the Central Limit Theorem for sequences of random variables which
are neither independent nor equally distributed. Its proof can be found in
\cite[Corollary 2.1.10]{duflo}.

\begin{theorem}
\label{teorema auxiliar}Let $\left\{  M_{n}\right\}  _{n\in\mathbb{N}}$ be a
real, square-integrable martingale adapted to a filtration $\left\{
\mathcal{F}_{n}\right\}  _{n\in\mathbb{N}}$. Let $\left\langle M\right\rangle
_{n}$ denote its angle bracket process. Let $\left\{  a_{n}\right\}
_{n\in\mathbb{N}}$ be an increasing sequence such that $a_{n}\rightarrow
\infty$ as $n\rightarrow\infty$.

\begin{description}
\item[(i)] If $\frac{\left\langle M\right\rangle _{n}}{a_{n}}\underset
{n\rightarrow\infty}{\longrightarrow}\sigma^{2}>0$ in probability, then
$\frac{M_{n}}{a_{n}}\underset{n\rightarrow\infty}{\longrightarrow}0$ a.s..

\item[(ii)] If Lindbergs's condition holds, that is,%
\[
\frac{1}{a_{n}}\sum_{k=1}^{n}\operatorname*{E}\left[  \left\vert M_{k}%
-M_{k-1}\right\vert ^{2}\mathbf{1}_{\{\left\vert M_{k}-M_{k-1}\right\vert
\geq\varepsilon\sqrt{a_{n}}\}}~|~\mathcal{F}_{k-1}\right]  \underset
{n\rightarrow\infty}{\longrightarrow}0\text{ in probability,}%
\]
then $\frac{M_{n}}{\sqrt{a_{n}}}\underset{n\rightarrow\infty}{\longrightarrow
}\mathcal{N}\left(  0,\sigma^{2}\right)  $ in law.\smallskip
\end{description}
\end{theorem}

\noindent\textbf{Acknowledgements.} The authors would like to thank Juan-Pablo
Ortega for his enlightening comments and suggestions. They also acknowledge
support from Centre de Recerca Matem\`{a}tica. Sebastian del Ba\~{n}o Rollin
is partially supported by Project MTM2006-01351, Ministerio de Educaci\'{o}n y Ciencia.

\end{document}